\documentclass[hidelinks,onefignum,onetabnum]{siamart250211}

\usepackage{setspace}
\usepackage{amsmath}
\usepackage{amssymb}
\usepackage{amsfonts}
\usepackage{mathrsfs}
\usepackage{bm}
\usepackage{bbm}
\usepackage[title]{appendix}
\usepackage{graphicx}
\usepackage{tikz}
\usepackage{pgfplots}
    \usetikzlibrary{spy}
\usepackage{caption}
\usepackage{subcaption}
\pgfplotsset{compat = 1.3}
\usepackage[margin=1in]{geometry}
\usepackage{comment}
\usepackage{multirow}
\usepackage{xcolor}
\usepackage{pgfplotstable, booktabs}
\usepackage{algorithm}
\usepackage{algpseudocode}
\usepackage{diagbox}
\usepackage{enumitem}
\usepackage{comment}
\usepackage{url}

\headers{Geometric $\alpha$SA multigrid for DG methods}{Y. Pan, M. Lindsey, P.-O. Persson}

\title{Geometric adaptive smoothed aggregation multigrid for discontinuous Galerkin discretisations\thanks{Submitted to the editors 17 April 2025.}}

\author{Yulong Pan\thanks{Department of Mathematics, University of California, and Mathematics Group, Lawrence Berkeley National Laboratory, CA, 94720, United States,
  (\email{yllpan@berkeley.edu}, \email{lindsey@berkeley.edu},
  \email{persson@berkeley.edu}).}
\and Michael Lindsey\footnotemark[2]
\and Per-Olof Persson\footnotemark[2]}

\begin{document}

\maketitle

\begin{abstract}
We present a geometric multigrid solver based on adaptive smoothed aggregation \cite{brezina2005adaptive} suitable for Discontinuous Galerkin (DG) discretisations. Mesh hierarchies are formed via domain decomposition techniques, and the method is applicable to fully unstructured meshes using arbitrary element shapes. Furthermore, the method can be employed for a wide range of commonly used DG numerical fluxes for first- and second-order PDEs including the Interior Penalty and the Local Discontinuous Galerkin methods. We demonstrate excellent and near mesh-independent convergence for a range of problems including the Poisson equation, and convection-diffusion for a range of P\'eclet numbers.
\end{abstract}



\section{Introduction}

Numerical methods for solving partial differential equations have become ubiquitous in both industry and academic settings for modeling a wide variety of physical phenomena of interest to modern practitioners. One such method which is the focus of this paper is the Discontinuous Galerkin (DG) method, which has become increasingly popular for simulations involving hyperbolic systems of conservation laws such as the Euler and Navier-Stokes equations. The popularity of the DG method is usually attributed to its high-order accuracy on arbitrarily unstructured meshes, in addition to its natural ability to stabilise numerical discretisations through the use of numerical fluxes.

A longstanding complaint about the DG method is its computational cost both in the discretisation and in the solver for resulting linear systems. Operators arising from DG discretisations are often more expensive to store than their counterparts from other methods and often require specialised solvers dependent on the specific choice of numerical flux. Moreover the solution of hyperbolic systems remains a significant challenge as solutions can contain artefacts from across the frequency spectrum.

One of the most important solver techniques currently employed is the multigrid method. Multigrid (MG) methods aim to iteratively solve a given linear equation via recursively combining degrees of freedom of the system into aggregates to obtain successively coarser representations of the problem. The main idea behind the multigrid algorithm is a divide-and-conquer principle, according to which different modes of the solution error are resolved separately at coarsened levels of representation of the original problem. The popularity of multigrid methods stems from their optimal linear asymptotic complexity when applied to elliptic operators, in particular the Laplacian.

Geometric multigrid (GMG) methods coarsen by directly combining neighbouring degrees of freedom on the underlying mesh, which is thus required as input to the solver in addition to the system matrix $A$. While the procedure for constructing GMG solvers on Cartesian meshes is well known, the extension to unstructured meshes has proven to be much more difficult. Furthermore coarsening strategies in GMG methods must be specialised for different element shapes. Various strategies have been employed over the years to tackle these difficulties, including remeshing \cite{adams1999parallel, chan1994domain} and element agglomeration \cite{dargaville2021comparison,ekstrom2010agglomeration,pan2022agglomeration}, but arguably it remain an open problem and is one main motivations for algebraic multigrid.

Algebraic multigrid (AMG) methods instead coarsen the problem using only the graph of the operator $A$. Degrees of freedom are combined in AMG according only to the non-zeros entries of $A$, bypassing the need for the mesh as input. Some popular variants of AMG include the classical Ruge-St{\"u}ben algorithm \cite{ruge1987algebraic}, smoothed aggregation \cite{vanek1996algebraic}, and adaptive AMG \cite{brezina2005adaptive,brezina2006adaptive}. While AMG solvers have been demonstrated to achieve good performance for numerous applications in Finite Elements
\cite{bochev2003improved,griebel2003algebraic,notay2016new}, its use in DG has been more limited so far. Recent work has aimed to extend the framework into DG \cite{antonietti2020algebraic,bastian2012algebraic,siefert2014algebraic,sivas2021air}, but many of these constructions are focused only on specific choices of numerical flux. 

AMG methods can be simpler to deploy due to their black-box nature, not requiring any information besides the system matrix $A$. However, such reliance on the structure of $A$ also means that AMG can be difficult to implement in a fully matrix-free fashion. One strategy recently introduced to alleviate this difficulty is low-order preconditioning, in which a geometric strategy is employed for a low-order grid overlaying a given high-order discretization on which AMG is applied \cite{chalmers2018low,heys2005algebraic,pazner2020efficient}. In this work we instead focus on building a geometric multigrid solver, although we adopt aspects of AMG methods to aid in generalising the method for unstructured mesh domains.

Though adaptive MG approaches have their origins~\cite{brezina2005adaptive,brezina2006adaptive} in the numerical PDE community, arguably their adoption has been most prominent in the lattice physics community, specifically lattice quantum chromodynamics (LQCD) \cite{Adaptive_MG_LQCD, MG_deflation_LQCD, alphaSA_LQCD}, as well as emerging applications in condensed matter physics \cite{MG_fermi} with similar mathematical structure. In these settings, the linear operator involves contributions that depend on highly disordered random vector or gauge fields, and functions in the near kernel space can look very non-smooth in the conventional sense.

We borrow heavily from adaptive Smoothed Aggregation ($\alpha$SA) AMG \cite{brezina2005adaptive} for the construction presented here. $\alpha$SA is a variant of Smoothed Aggregation (SA) AMG in which the construction of the multigrid hierarchy is automated through an adaptive process that tailors the solver to different problems without the need for extensive manual tuning by the user. The generality of the $\alpha$SA framework allows our construction to be used for a wide range of DG numerical fluxes, without specific tuning to each case.

In addition, we utilise an SVD step in the adaptive process for constructing the hierarchy, as well as an adaptive smoothing algorithm in the multigrid solver. These extensions to the original $\alpha$SA algorithm were used in previous work \cite{chow2006aggregation, d2019improving} to improve the performance of resulting the adaptive AMG solves. We leverage these constructions to achieve consistent consistent performance across different problems, mesh discretizations, and choices of numerical DG flux without requiring separate parameter hand-tuning for different cases.

A major point distinguishing our method from this previous work is that ours is ultimately a GMG solver, in that hierarchies are constructed using mesh information rather than using the adjacency graph of the system matrix $A$. As such, mesh information (i.e., specifying which indices correspond to which mesh elements) must be supplied as an input by the user. This avoids the need for parameter tuning in algebraic coarsening as well and allows the method to be implemented naturally in a matrix-free manner, rendering it suitable for use on distributed computing architectures. 

To construct mesh hierarchies we can employ domain decomposition methods popular in Finite Element methods to form element aggregates. Specifically we use the popular METIS package \cite{karypis1997metis} in this work. In principle METIS can be applied using non-uniform edge weights~\cite{karypis1997metis}, and in future work we intend to use this capacity to tackle problems with significant mesh anisotropies. 

We furthermore leverage knowledge of the underlying mesh to enable coarsening within individual elements at each level of the multigrid hierarchy, akin to p-multigrid methods employed for DG \cite{fidkowski2005p,persson2008newton}. We refer to such intra-element coarsening as an h*-multigrid step, as opposed to standard inter-element h-multigrid coarsening strategies usually applied in geometric MG solvers. The use of this additional coarsening strategy allows for the construction of efficient solvers for DG discretisations at higher p-refinements.

We apply our method both as a iterative solver and as a preconditioner to the Conjugate Gradient/Generalised Minimum Residual methods. The method is employed on unstructured hexahedral and tetrahedral meshes, and we find the performance on both element types to be extremely comparable without modifications in the solver construction. We test the method on convection-diffusion problems at a range of P\'eclet numbers and we show that the number of iterations to convergence using our solver remains largely constant across the board, in contrast to classical GMG and AMG techniques which tend to be specialised only for diffusion-dominated problems.

The paper is structured as follows. First we provide a brief review of the DG method in Section \ref{sect:DG}, as well as a review of the relevant details of smoothed aggregation and adaptive smoothed aggregation AMG in Section \ref{sect:AMG}. Details of our method are then outlined in Section \ref{sect:GMG}. Numerical examples employing the method on a range of test problems are shown in Section \ref{sect:examples} before we finally conclude in Section \ref{sect:conclusion}.

\section{Discontinuous Galerkin formulation} \label{sect:DG}
We review the details of the Discontinuous Galerkin (DG) method for a general scalar equation
\begin{equation}
    \nabla \cdot D(u,\nabla u) = f
\end{equation}
in some domain $\Omega \in \mathbb{R}^d$ where $d$ denotes the dimension and $D$ some function of $u$ and its gradient $\nabla u$. To begin, the domain $\Omega$ is discretised using a mesh $\mathcal{T}_h = \{ K_i \}$ consisting of distinct elements $K_i$ such that $K_i \cap K_j = \emptyset$ and $\cup_i K_i = \Omega$. A broken function space is defined on the set of elements as
\begin{equation}
    V(\mathcal{T}_h) = \{ v \in L^2(\Omega) : v|_{K_i} \in H^1(K_i),~ \forall K_i \in \mathcal{T}_h \}
\end{equation}
A suitable Finite Element space $V_h \subset V (\mathcal{T}_h) $ is then chosen by the user based on a further projection of $H^1(K_i)$ on each element to a finite-dimensional space. A popular choice for $V_h$ is the broken polynomial space of degree of at most $p$.

To obtain an approximate solution $u \in V_h$, a weak form of the problem is obtained by applying a standard Galerkin procedure to the equation of interest
\begin{equation}
    \sum_i \int_{K_i} [ v \nabla \cdot D (u, \nabla u) - vf ] ~dx = 0, 
\end{equation}
where $v \in V_h$ is some test function. For a general DG method, integration by parts is typically applied at this point, yielding 
\begin{equation}
\label{eq:ibp}
    \sum_i \int_{\partial K_i} v D(u, \nabla u) \cdot \bm{n} ~ds = \sum_i \int_{K_i} [\nabla v \cdot  D  (u, \nabla u)  + vf ]~dx
\end{equation}
where $\bm{n}$ denotes the outward pointing normal vector, and the hat-notation has been introduced to indicate the numerical flux on an element boundary. The numerical flux is needed to ensure that the problem is well defined on element boundaries, where $u \in V_h$ is in general multi-valued. 

The choice of numerical flux depends on the exact problem considered: for first order systems some common choices include the Godunov and Roe fluxes, whilst for second order problems popular choices include the Interior Penalty \cite{arnold1982interior} and Local Discontinuous Galerkin \cite{cockburn1998local} fluxes. In this work we do not make any assumption about the choice of flux, developing our solver for general DG discretisations of both first- and second-order problems, in which the specification of the numerical flux can be arbitrary.

In the case where $D$ is linear, \eqref{eq:ibp} ultimately yields a linear system of the form
\begin{equation} \label{eq:Auf}
    Au = f
\end{equation}
where we have chosen to absorb the mass matrix into the right hand side vector $f$. For a nonlinear $D$, one typically uses Newton's method and solves a sequence of such linear systems. Ultimately the key computational difficulty is the solution of the linear system involving an operator $A$ of size $n \times n$.

\section{Algebraic multigrid methods} \label{sect:AMG}
In this section we provide a brief overview of some previous work on algebraic multigrid (AMG) methods relevant for our discussions.

\subsection{General algorithm}
In algebraic multigrid (AMG), node aggregates are formed to coarsen the linear problem $Au=f$ using the graph of the operator $A$ without any knowledge of the underlying mesh. While we do not delve fully into the details of algebraic aggregation the general idea of AMG is to infer, from the weighted graph induced by the non-zero entries of $A$, local clusters of degrees of freedom which are strongly connected to one another. These neighbourhoods are then combined into an aggregate and the aggregates collected into a level. On each level $k$, denoted $L_k$, a new operator $A_k$ is then formed and the process repeated recursively until the operator on the final level $A_N$ is sufficiently small to be inverted directly.

To define each of these levels, the following operators are constructed:
\begin{enumerate}
    \item Between two successive levels $k$ and $k+1$, a prolongation operator $T_{k+1}^{k}$ projecting vectors defined on level $k+1$ to level $k$.
    \item Between two successive levels $k$ and $k+1$, a restriction operator $R_{k}^{k+1}$ projecting vectors defined on level $k$ to level $k+1$.
    \item On each level $k$, the reduced operator $A_k$ given by $A_{k+1} = R_{k}^{k+1} A_{k} T^k_{k+1}$, where $A_0 = A$.
    \item A smoothing operator $\tilde{A}_k$ that approximates $A_k$ but is much easier to invert than $A_k$ itself.
\end{enumerate}
The exact choices of prolongation, restriction, and smoothing operators are abstracted away for the moment and are expanded upon in the following sections. The smoother $\tilde{A}_k$ is used at level $L_k$ to approximately solve the equation $A_k u = b_k$ via an iterative relaxation procedure (Algorithm \ref{alg:smooth}).
\begin{algorithm}[h]
\caption{Smooth $(u_0, A_k, \tilde{A}_k, b_k, \text{stop($\cdot$)})$}
\begin{algorithmic}[1]
    \State Define $x^{(0)} = u_0, ~m=0$
    \While {!\,stop($\cdot$)}
        \State Set $m = m+1$
        \State Set $x^{(m)} = x^{m-1} + \tilde{A}_k^{-1} ( b_k - A_kx^{(m-1)})$
    \EndWhile
    \State Return $x^{(m)}$
\end{algorithmic}
\label{alg:smooth}
\end{algorithm}

In Algorithm \ref{alg:smooth}, $u_0$ denotes the initial guess to the iterative smoothing procedure. The stopping criterion stop$(\cdot)$ can be tuned for specific use cases and is specified in more detail in subsequent algorithms where the smoothing operation is applied.

\begin{algorithm}[h]
\caption{Multigrid V-Cycle $(u_0, A, f)$}\label{alg:mg}
\begin{algorithmic}[1]
    \State Define $b_0 = f - Au_0$, where $u_0$ denotes the initial guess
    \For {$k=0$ to $N-1$}
        \State Set $u_k$ = Smooth($0, A_k, \tilde{A}_k, b_k, m \geq t_\text{pre})$ \Comment{pre-smoothing}
        \State Calculate residual $r_k = b_k - A_k u_k$
        \State Define $b_{k+1} = R_{k+1}^k r_k$ \Comment{restrict residual}
    \EndFor
    \State Solve $u_N = A_{N}^{-1} b_N$ \Comment{bottom level solve}
    \For {$k = N-1$ downto $0$}
        \State Update approximate solution $u_k = u_k + T_{k}^{k+1} u_{k+1}$ \Comment{prolong solution}
        \State Set $u_k$ = Smooth($u_k, A_k, \tilde{A}_k, r_k - Au_k, m \geq t_\text{post}$) \Comment{post-smoothing}
    \EndFor
    \State If $\| b - A u_0 \| > \epsilon$, go back to Line 2
    \State Return $u_0$
\end{algorithmic}
\label{alg:MultigridVCycle}
\end{algorithm}

With the smoothing operation defined, AMG typically solves the equation $Au=b$ using the V-Cycle meta-algorithm specified in Algorithm~\ref{alg:MultigridVCycle}, which is also used by geometric multigrid methods.
In the smoothing steps specified on Lines 3 and 10, the stopping criterion is specified such that exactly $t_\text{pre}$ and $t_\text{post}$ smoothing iterations are applied, respectively. This algorithm is known as the multigrid V-cycle as the algorithm first traverses down the hierarchy before retracing its steps back up level-by-level, applying the smoother before each descent step and after each ascent step (`pre-smoothing' and `post-smoothing,' respectively).

\subsection{Smoothed aggregation}
Smoothed aggregation (SA) multigrid \cite{vanek1996algebraic} is a variant of AMG which constructs the restriction operator given a set of vectors defining the coarsening on each level of the hierarchy. At any given level $k$, standard AMG techniques are first applied to the level operators $A_k$ to form `aggregates' $\mathcal{A}_i^k$ denoting a partition of the degrees of freedom at the level $k$ into non-overlapping subsets. The SA algorithm then takes in a matrix $B^{k+1}$ whose columns represent the coarse space of the following level $k+1$ and defines in terms of $B^{k+1}$ and the aggregates
\begin{equation} \label{eq:splitB}
    P_{k+1}^{k} = 
    \begin{pmatrix}
    \mathrm{orth} ( B^{k+1}[\mathcal{A}_{1}^{k},:] ) & &  \\
    & \mathrm{orth} ( B^{k+1}[\mathcal{A}_{2}^{k},:] ) & & \\
    & & \ddots
    \end{pmatrix}.
\end{equation}
Here $B^{k+1}[\mathcal{A}_{i}^{k},:]$ denotes the rows of $B^{k+1}$ corresponding to the indices of degrees of freedom in the aggregate $\mathcal{A}_{i}^{k}$. This operation carves up $B^{k+1}$ into blocks of rows which are then orthogonalised and unfolded into a block-diagonal matrix which by construction is an isometry, i.e., satisfies $(P_{k+1}^k)^\top P_{k+1}^k = I$. In practice the orthogonalization is achieved by taking the Q factor of a QR factorisation.

As input, SA only requires the matrix $B^1$ at the top level. Subsequent $B^{k+1}$ are then formed as
\begin{equation} \label{eq:SA_Bk}
    B^{k+1} := (P_{k+1}^k)^\top B^k.
\end{equation}

The operator $P_{k+1}^k$ is related to the prolongation operator $T_{k+1}^k$, which is constructed in the SA method as 
\begin{equation} \label{eq:SA_prolong}
    T^k_{k+1} = S_{k} P_{k+1}^k
\end{equation}
where $S_{k}$ is a smoothing operator that smooths out the jumps across aggregates introduced by the splitting procedure~\eqref{eq:splitB}. A simple choice for $S_{k}$ is given by 
\begin{equation}
\label{eq:SA_Sk}
    S_{k} = I - \tilde{A}_k^{-1} A_k
\end{equation}
which is the same smoother used in the multigrid V-cycle. With the prolongation operator defined, the restriction operator $R^{k+1}_k$ for SA is related to the prolongation simply by
\begin{equation}
\label{eq:SA_restrict}
    R^{k+1}_k = (T_{k+1}^k)^\top.
\end{equation}

\subsection{Adaptive smoothed aggregation}
A key difficulty faced by the smoothed aggregation algorithm is the specification of the candidate matrix $B^1$ which must be specified for each problem. The main idea behind adaptive smoothed aggregation ($\alpha$SA) AMG \cite{brezina2005adaptive} is to automate this process by using the multigrid smoother $\tilde{A}$ to find a suitable candidate matrix. 

The philosophy behind the method is to use the smoother to identify a near kernel space that can be used to define subsequent levels of the multigrid hierarchy. At the top level, a random $n \times r$ matrix $B^1 = (b_1^1~ ... ~b_r^1)$ with independent standard gaussian entries defines an initial guess for $B^1$, the columns of which are refined via applications of the smoother $\tilde{A}$ to approximately solve the homogeneous equation $A x = 0$. The resulting columns $B^1 = (b_1^1, \ldots, b_r^1)$ are assembled into the matrix $B^1$, which is then propagated to subsequent levels via \eqref{eq:SA_Bk}. These can then be adaptively refined by applying the smoother $\tilde{A}_k$ at each level. 

In general the idea of $\alpha$SA is that after applying the smoother for sufficiently many steps, the problem can be projected to the degrees of freedom that are difficult for the smoother $\tilde{A}_k$ to resolve effectively. In terms of $B^k$, the prolongation and restriction operators are defined just as in SA, following \eqref{eq:SA_prolong} and \eqref{eq:SA_restrict}. The resulting algorithm is summarized in Algorithm~\ref{alg:adaptiveSA}.

\begin{algorithm}[h]
\caption{Adaptive smoothed aggregation $(A)$}\label{alg:adaptiveSA}
\begin{algorithmic}[1]
    \State Set $A_0 = A$
    \State Choose $B^1 = (b_1^1 ~...~ b_r^1)$ to be a random gaussian $n \times r$ matrix
    \For {$k = 1$ to $N-1$ \label{step:loop}} 
        \State Construct smoother $\tilde{A}_k$ from $A_k$
        \For {$j = 1$ to $r$}
            \State Update $b_j^k$ $\leftarrow$ Smooth($b_j^k, A_k, \tilde{A}_k, 0, m \geq t_\text{init}$)
        \EndFor
        \State Assemble $B^k = (b_1^k, \ldots, b_r^k)$
        \State Construct aggregates $\{ \mathcal{A}^k_i \}$ in level $L_k$ using $A_{k}$
        \State Use $B^k$ to form $P_{k+1}^{k}$ via Equation (\ref{eq:splitB})
        \State Construct prolongation operator $T_{k+1}^{k} = S_{k} P_{k+1}^{k}$ and restriction operator $R_{k}^{k+1} = (T_{k+1}^{k})^\top$
        \State Construct next level operator $A_{k+1} = R_{k}^{k+1} A_{k} T_{k+1}^{k}$ 
        \State Initialize $B^{k+1}$ using Equation (\ref{eq:SA_Bk})
    \EndFor
    \State Construct multigrid solver using prolongation and restriction operators $T_{k+1}^k$ and $R_k^{k+1}$, $k=0,\ldots,N-1$
    \State Pick a random vector $\tilde{x}$
    \State Use multigrid solver on $Ax=0$ with initial guess $\tilde{x}$, i.e., run $\text{Multigrid}(\tilde{x},A,0)$
    \State If performance of multigrid solve is unsatisfactory, go back to Line \ref{step:loop}
\end{algorithmic}
\end{algorithm}

To evaluate the hierarchy, the AMG solver constructed using $\alpha$SA is tested on the homogeneous problem $Ax = 0$, with a random initial guess $\tilde{x}$. The procedure can be refined until its performance on the homogeneous performance is satisfactory. As the algorithm makes no prior assumptions about the structure of $A$, the hope is that this adaptive procedure can be applicable to a wide range of problems without the need for extensive fine-tuning.

\section{Geometric adaptive smoothed aggregation multigrid for DG} \label{sect:GMG}
We borrow from the framework of adaptive smoothed aggregation AMG to construct a geometric multigrid solver for Discontinuous Galerkin (DG) discretisations. In contrast with AMG, node aggregations are formed using the mesh, not the adjacency graph structure of the linear operator $A$. As such this method can be implemented without knowledge of the sparsity pattern of $A$.

\subsection{Geometric hierarchies}
Given a DG mesh $\mathcal{T}_h$ for a domain contained in $\mathbb{R}^d$, we construct an associated unweighted graph $\mathcal{G}$ where vertices correspond to the elements $K_i$. Two vertices are connected by an edge if and only if the corresponding elements $K_i, K_j$ share a face.

We define the aggregates $\mathcal{A}_i^k$ forming each level $k$ of the hierarchy by applying a hierarchical clustering algorithm to this graph. 
For this work we use the popular package METIS \cite{karypis1997metis} to perform this task. Other techniques that have been used in the literature for this task include Graph Neural Networks \cite{antonietti2024agglomeration} and greedy element clustering \cite{pan2022agglomeration}. We do not undertake a comparative study of different agglomeration strategies as we have found that the performance of METIS to be highly effective for the examples considered in this work.

Specifically, at the coarsest bottom level which we denote as the $N^\text{th}$ level, we cluster the mesh graph into $2^d$ connected components such that each component consists of roughly the same number of elements. These components define the aggregates $\mathcal{A}_i^N$ forming the bottom level. Proceeding inductively, the aggregates at each level are each split into $2^d$ connected components to define the aggregates at the subsequently higher level. If an aggregate already consists of fewer than $2^d$ elements, then it is simply split into singletons. The procedure is repeated until each aggregate consists of a single element, defining the top level $0$ of the hierarchy.



\begin{figure}
    \centering
    \includegraphics[scale=0.4]{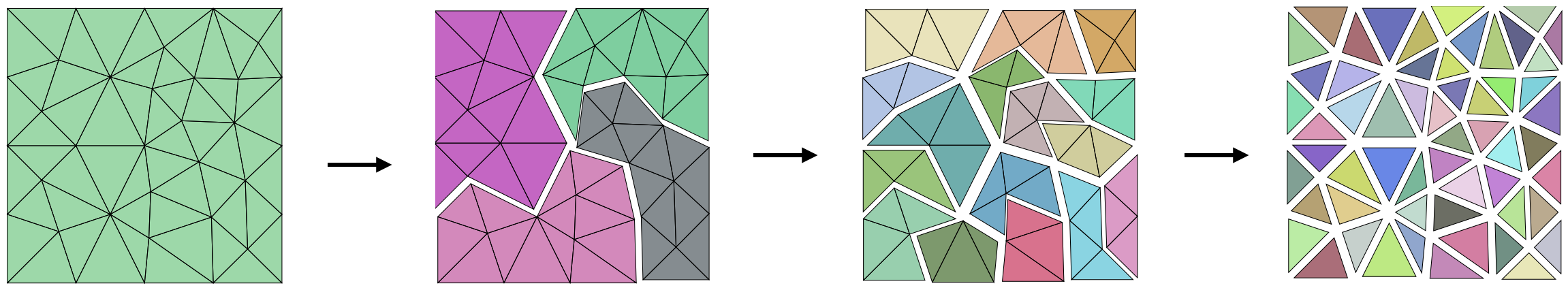}
    \caption{Schematic of splitting procedure to form aggregates. Levels are formed in reverse where mesh is first split into four blocks, then each is further split into four blocks. The finest level $L_0$ is defined by the last collection of blocks where each aggregate consists of only one element.}
    \label{fig:metis_split}
\end{figure}

A schematic of the hierarchy construction is shown in Figure \ref{fig:metis_split} for a rectangular mesh, though note that the procedure can be applied as described to unstructured meshes.

\subsection{Prolongation/restriction operators}
With the geometric hierarchy formed, we turn to the construction of the prolongation and restriction operators. At each $k$-th level following the $\alpha$SA framework we take a candidate matrix $B^k = (b_1^k ~...~ b_{r}^k)$ of size $n \times r$, whose columns we smooth according to the homogeneous equation $\tilde{A}_k x = 0$:
\begin{equation} \label{eq:bsmooth}
    b_j^k \leftarrow \text{Smooth}\bigg( b_j^k, A_k, \tilde{A}_k,0, \frac{\| x^{(m)} \|}{\| x^{(m)+1} \|} < 1 + \gamma \bigg)
\end{equation}
The stopping criterion checks whether the relative change in norm between successive iterates in the smoothing operation falls below a certain tolerance, with the idea being that once convergence has stagnated, each of the columns of $B^k$ should consist of mainly the modes that are not well resolved by the smoother. In our implementation the tolerance is set to $\gamma=0.03$, although we have found in our testing for the performance to not be overly sensitive to the specific choice of $\gamma$. Further details about the smoother will be offered in Section \ref{sect:smoothing}.

We follow a similar methodology to \cite{chow2006aggregation,d2019improving} where we over-sample the number of columns $r$ in $B^k$ before trimming them down using a singular value decomposition, as we elaborate upon here. At the top level, we choose the number of columns $r$ to equal the median number of degrees of freedom in a single aggregate at the first level. The matrix $B^1$ is then formed as a random Gaussian matrix of $n \times r$ entries, and its columns smoothed according to Equation \eqref{eq:bsmooth}. For lower levels $k$, the candidate matrix is initialised using the restriction operator according to Equation \eqref{eq:SA_Bk}, then smoothed as before. 

Then we form the coarsening operator in a similar fashion to $\alpha$SA. For each aggregate $\mathcal{A}_i^{k+1}$, we first isolate the rows of $B^1$ corresponding to the degrees of freedom of all elements in the aggregate. However instead of orthogonalizing all of the columns of $B^k [\mathcal{A}_i^{k+1}, \, : \,]$ to define a suitable block of $P_{k+1}^k$ (cf. \eqref{eq:splitB} for comparison), we perform a singular value decomposition and pick a subset of the left singular vectors to span the coarsened space for the aggregate. We denote the number of vectors kept at each level by $s_i^k$. We elaborate on the choice of this value in the following Section \ref{sect:hp}.

In summary, $B^k$ is used to form $P_{k+1}^k$ by the following procedure:
\begin{equation} \label{eq:splitB_dg}
    B^{k}
    \rightarrow
    \begin{pmatrix}
    B^{k}[\mathcal{A}_{1}^{k+1},\,:\,] = U^k_{1}\Sigma^k_{1} (V^k_{1})^\top & &  \\
    & \ddots
    \end{pmatrix}
    \rightarrow
    P_{1}^k =
    \begin{pmatrix}
    U^k_{1}[\,: \, ,1\,\mathrm{:}\,s_1^k] & &  \\
    & \ddots
    \end{pmatrix}.
\end{equation}
Meanwhile, $T_{k+1}^k$ and $R_k^{k+1}$ are determined in terms of $P_{k+1}^k$ following \eqref{eq:SA_prolong} and \eqref{eq:SA_restrict} as in SA. For the choice of smoothing operator $S_1$ used to define the prolongation operator, cf. \eqref{eq:SA_Sk}, again we use the same smoother $\tilde{A}_k$ as in the multigrid V-cycle. Concretely, then, we have:
\begin{equation}
\label{eq:TR0}
    T^k_{1} = (I-\tilde{A}^{-1}_k A_k) P_{k+1}^k, \quad R_k^{k+1} = (T^k_{k+1})^\top.
\end{equation}
Moreover the operator at the next level of the hierarchy is defined as
\begin{equation}
    A_{k+1} = R_k^{k+1} A_k T_{k+1}^k,
\end{equation}
together with an associated $\tilde{A}_1$.

Note that the construction of $P_{k+1}^k$ can be viewed as a randomized sketching approach for computing the SVD~\cite{Martinsson_Tropp_2020, murray2023randomizednumericallinearalgebra}. To wit, if $F$ denotes the function such that the operator $F(A)$ implements the action of the smoother, then at the top level, for each $i$, we are precisely performing a randomized SVD on the matrix $F(A) [\mathcal{A}_i^{k+1}, \,:\,]$, since in fact $ B^k [\mathcal{A}_i^{k+1}, \,:\,] = ( F(A) [\mathcal{A}_i^k, \,:\,] ) \,  Z$, where $Z$ is a random gaussian sketch matrix. Consistent with the perspective from randomized numerical linear algebra, even if the target truncation rank is $s_i^k$, it makes sense to consider a sketch matrix with $r > s_i^k$ columns in order to determine the optimal truncation, over-sampling by a constant (or logarithmically growing) factor.

\subsection{Intra- vs inter-aggregate coarsening} \label{sect:hp}
While in the original $\alpha$SA methodology coarsening is only performed across aggregates, due to the block structure of the algorithm we may also consider additional coarsening steps within each aggregate $\mathcal{A}_i^{k}$. This is analogous to the methodology used by p-multigrid solvers. However we note that such intermediate steps are optional.

\begin{figure}
    \centering
    \includegraphics[scale=0.35]{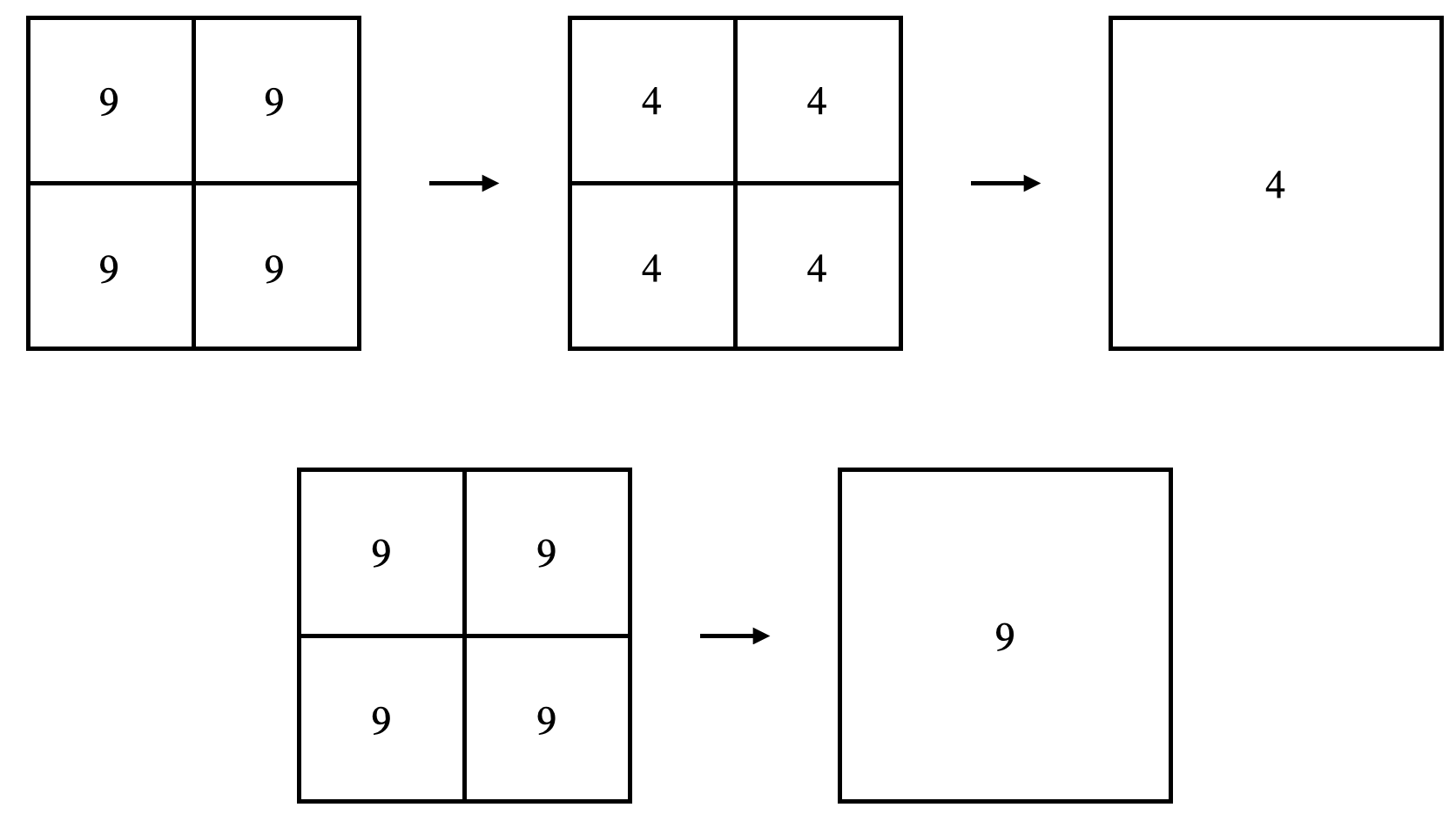}
    \caption{Comparison of h*- and h- multigrid approaches. Numbers correspond to degree of freedom count in each element. On the top row, the h*-multigrid approach is shown: intra-element coarsening is first applied in each element prior to geometric coarsening. On the bottom row, the h-multigrid approach is shown: geometric coarsening is applied directly without pre-coarsening within each element. 
    }
    \label{fig:hp_h_mg}
\end{figure}

In this work we refer to the first case where coarsening is performed across aggregates as an h-multigrid step, and the latter case where coarsening is applied in an intra-aggregate manner as an h*-multigrid step. A schematic comparing both the h*- and h-multigrid approaches is shown in Figure \ref{fig:hp_h_mg}. The main function of including the intra-aggegate coarsening steps is to cut down the number of degrees of freedom within an aggregate prior to geometric coarsening. In our numerical examples, we demonstrate that such inclusions can yield a faster reduction of the operator complexity in the multigrid hierarchy.

Within our formalism, in order to introduce a h*-multigrid step at the $k$-th level, we simply duplicate the level and its aggregates in the multigrid hierarchy. Then according to Equation \eqref{eq:splitB_dg}, the spanning set for each coarsened space will be formed only using intra-aggregate degrees of freedom. 

We note that many p-multigrid coarsening strategies have been employed in the literature. In \cite{fidkowski2005p}, degrees of freedom are eliminated by decreasing the polynomial order in each element by one at each multigrid level whilst in \cite{persson2008newton} the polynomial order is reduced more aggressively by roughly by a factor of two at each level. By contrast, in the spirit of adaptive AMG techniques, we will not coarsen explicitly according to polynomial degree in the h*-multigrid approach, but rather we will reduce the number of degrees of freedom by choosing the coarsened bases in each element adaptively.

Specifically, in \emph{either} an h- or h*-multigrid step, we look to decrease the number of degrees of freedom in each element according to the singular values of each block in Equation \eqref{eq:splitB_dg}
\begin{equation}
    B^{k}[\mathcal{A}_{j}^{k+1},:] = U^k_{j} \Sigma^k_{j} (V^k_{j})^\top
    \rightarrow 
    U^k_{j}[\,: \, ,1\,\mathrm{:}\,s_j^k]
\end{equation}
where the number of left singular vectors kept $s_j^k$ is given by
\begin{equation} \label{eq:sik}
    s_j^k = \begin{cases}
        |\mathcal{A}^k_j| \, / \, n_\text{cut}, &\text{if } \mathcal{A}_j^k \neq \mathcal{A}_j^{k+1} ~(\text{h-multigrid}), \\
        \left|\{ \sigma_l \in \text{diag}(\Sigma_j^k) : \sigma_l > \delta \sigma_1 \} \right|, &\text{if } \mathcal{A}_j^k = \mathcal{A}_j^{k+1} ~(\text{h*-multigrid}).
    \end{cases}
\end{equation}

In an h-multigrid step, we look to cut the number of functions in each aggregate, denoted $|\mathcal{A}^k_j|$, by a constant factor $n_\text{cut}$ such that the rate at which degrees of freedom within each aggregate decrease roughly matches the number of aggregates combined in a geometric coarsening step. This ensures that the size of $A^{k+1}$ decreases roughly by a constant factor with each h-multigrid step, in line with traditional h-multigrid methods. Although ideally we would have $n_\text{cut}$ match exactly the $2^d$ aggregates combined at each level of the mesh coarsening procedure, we observe in numerical tests in Section \ref{sec:ncut} that we need to set $n_\text{cut}$ slightly smaller than this in order to achieve optimal multigrid performance. As such, we set a default value of $n_\text{cut} = 2^d - (d-1)$ in our implementation, which we find experimentally to achieve good solver performance. However, we expect that the optimal value of this parameter may depend on the specific application.

Meanwhile, in an h*-multigrid step, we look to cut the number of functions according to singular values of suitable blocks of $B^{k}$. In this case we do not have the same geometric considerations as in the case of an h-multigrid step. Specifically, we choose to keep left singular vectors of $B^{k}[\mathcal{A}_{j}^{k+1},:]$ corresponding to singular values within a threshold factor $\delta$ of the largest singular value $\sigma_1$ of the block. For our implementation, we set $\delta = 10^{-3}$.

Other than the choice of aggregates and the value of $s^k_j$, there is no algorithmic difference between the intra-element and geometric coarsening levels. We emphasise again that the h*-multigrid steps are optional; we consider both h*- and h-multigrid approaches in our numerical experiments in Section \ref{sect:examples}. In this work, when h*-multigrid is employed, we only apply it at the top level of the hierarchy, similar to the hp-multigrid methodology. While more h*-multigrid steps could be added on subsequent levels, we found in practice that our choice was able to achieve good solver performance without requiring any additional levels in the hierarchy.

\subsection{Adaptive smoothing} \label{sect:smoothing}
As explained above, the multigrid construction requires us to specify a choice of $\tilde{A}_k$ given the operator $A_k$ at each level of the hierarchy. For this work we focus on choices of the form 
\begin{equation} \label{eq:smoother_choice}
    \tilde{A}_k = \frac{1}{\omega_k} P(A_k), 
\end{equation}
where $\omega_k$ is a damping factor to be specified below and $P(A_k)$ is some matrix function of $A_k$. Some common choices of $P$ for DG methods include:
\begin{enumerate}
    \item Richardson: $P(A_k) = I$
    \item Block Jacobi: $P(A_k) = D_k$, the block diagonals of $A_k$
    \item Block Gauss-Seidel: $P(A_k) = L_k$, the lower block triangle of $A_k$
\end{enumerate}
For classical multigrid methods, the damping factor is chosen such that the spectrum of the matrix
\begin{equation} \label{eq:smoother_op}
    S_k = I - \omega_k P(A_k)^{-1} A_k,
\end{equation}
lies roughly in the range $[-\frac{1}{3},1]$ with eigenvectors corresponding to high-frequency modes of $A_k$ having eigenvalues contained roughly within the interval $[-\frac{1}{3},\frac{1}{3}]$. For the model problem of a Laplacian on a Cartesian grid, the optimal choice of damping factors achieving these desiderata are well-known for the block Jacobi and block Gauss-Seidel methods: $\omega = \frac{2}{3}$ and $\omega=1$ respectively. These values are often employed in practice on more complicated problems although it is less clear whether they are truly optimal in those cases.

We employ the simple adaptive procedure to determine the damping factor for each given $A_k$ and $P(A_k)$ commonly used in both SA and $\alpha$SA type implementations. We proceed by estimating the spectral radius of the operator $P(A_k)^{-1} A_k$ appearing in \eqref{eq:smoother_op} via the power iteration in Algorithm \ref{alg:power}. In our implementation we set the number of iterations to $q=3$.

\begin{algorithm}[h]
\caption{Estimate spectral norm $(A_k, P(A_k))$}
\begin{algorithmic}[1]
    \State Pick a random vector $x$ of unit norm
    \For {$i=1$ to $q$}
        \State Set $x = Ax$
        \State Set $x = P(A_k)^{-1}x$
        \State Normalise $x = x / \| x \|$
    \EndFor
    \State Return $\| P(A_k)^{-1} A_k x \|$
\end{algorithmic}
\label{alg:power}
\end{algorithm}

Given an estimate to the spectral radius $\tilde{\rho}_k \approx \rho(P(A_k)^{-1} A_k)$ we then set the damping factor $\omega_k$ to be
\begin{equation}
\label{eq:omega}
    \omega_k = \frac{4}{3} \cdot \frac{1}{\tilde{\rho}_k}
\end{equation}
and the corresponding smoother $\tilde{A}_k$ according to \eqref{eq:smoother_choice}. This choice aims to mimic smoothers for the standard model Laplacian as closely as possible using a damped block Jacobi smoother where $P(A_k) = D_k$. We demonstrate through our numerical tests that it works well in practice more generally.

\subsection{Algorithm summary}
We summarise the algorithm here to complete our discussion of the method. As input to construct the multigrid hierarchy and its operators we require the following information:
\begin{enumerate}
    \item $A$, the system operator of size $n \times n$.
    \item A graph $\mathcal{G}$ of the underlying mesh, where nodes correspond to elements and edges connect neighbouring elements.
\end{enumerate}
With these specified by the user the algorithm proceeds as outlined in Algorithm \ref{alg:full_mg}. Once all the operators in the hierarchy have been formed, the standard multigrid V-cycle in Algorithm \ref{alg:mg} can be applied without any modifications to solve the linear system $Au=f$.

\begin{algorithm}[!htbp]
\caption{Geometric $\alpha$SA multigrid $(A, \mathcal{G})$}\label{alg:full_mg}
\begin{algorithmic}[1]
    \State \emph{\# Geometric hierarchy}
    \State Form levels $\{ L_0, \ldots,L_N \}$ with aggregates $\{\mathcal{A}_i^k\}$ by applying METIS to the mesh graph $\mathcal{G}$
    \State Set $A_0 = A$
    \State
    \State \emph{\# Intra-element coarsening}
    \For {$k = 1$ to $N-1$}
    \If {h*-multigrid at level $L_k$}
        \State Replicate level $L_k$, i.e.,
        \begin{itemize}[leftmargin=.5in]
            \item[--] Relabel levels $\{ ...,L_k,L_{k+1},L_{k+2},... \} \leftarrow \{ ...,L_k,L_{k},L_{k+1},... \}$ 
            \item[--] Set $N \leftarrow N + 1$
        \end{itemize}
    \EndIf
    \EndFor
    \State
    \State \emph{\# Prolongation/restriction}
    \State Construct $B^1 = (b_1^1 ~...~ b_r^1)$ as a random gaussian $n \times r$ matrix
    \For {$k = 1$ to $N-1$}
    \State Estimate $\tilde{\rho}_k$ by applying Algorithm \ref{alg:power} to $(A_k, P(A_k))$
    \State Define $\omega_k$ using $\tilde{\rho}_k$ and Equation \eqref{eq:omega}
    \State Construct $\tilde{A}_k$ using $\omega_k$ and Equation \eqref{eq:smoother_choice}
    \For {$j = 1$ to $r$}
        \State Update $b_j^{k}$ $\leftarrow$ Smooth($b_j^{k}, A_k, \tilde{A}_k, 0, \frac{\| x^{(m)} \|}{\| x^{(m+1)} \|} < 1 + \gamma$)
    \EndFor
    \State Use $B^{k}$ to form $P_{k+1}^{k}$ via Equations \eqref{eq:splitB_dg} and \eqref{eq:sik}
    \State Construct prolongation operator $T_1^0 = (I - \tilde{A}_k^{-1} A_k) P_{k+1}^k$ and restriction operator $R_k^{k+1} = (T_{k+1}^k)^\top$
    \State Construct next level operator $A_{k} = R^{k+1}_k A_k T_{k+1}^k$ 
    \State Initialise $B^{k+1} = R_k^{k+1} B^{k} $
    \EndFor

    \State
    \State Return operators $A_k, \tilde{A}_k$ for all $k=0,\ldots,N$ and $T_{k+1}^{k}, R_{k}^{k+1}$ for all $k=0,\ldots,N-1$
\end{algorithmic}
\end{algorithm}

\subsection{Matrix-free implementation} \label{sect:matrix_free}
In its current formulation, while the adaptive multigrid construction guarantees that operators at successive levels $A_k$ ($k=0,1,...,N$) decrease monotonically in size, the same is not necessarily true of the number of non-zero entries. In numerical experiments presented below in Section \ref{sect:examples}, we demonstrate that for lower values of $n_\text{cut}$ in Equation \eqref{eq:sik}, significant increases in the number of non-zero entries can in fact be observed at higher levels in the multigrid hierarchy.

To counteract this we implement the multigrid algorithm such that the operators $A_k$ at the first $\kappa \in \{0,\ldots, N \}$ levels are not stored explicitly. Instead we can construct the multi-level restriction and prolongation operators between two levels $k_1 < k_2$ as
\begin{equation}
\label{eq:matfree1}
    T_{k_2}^{k_1} = \prod_{j=k_1}^{k_2} T_{j+1}^j, ~R_{k_1}^{k_2} =\bigg( T_{k_2}^{k_1} \bigg)^T
\end{equation}
and the operator at level $k$ as
\begin{equation}
\label{eq:matfree2}
    A_k = R^k_0 A_0 T_k^0.
\end{equation}
The multi-level transfer operators can be either stored explicitly or unrolled specifically depending on the memory restrictions of the user. Then matrix-vector multiplications by $A_k$ can be performed by applying these multi-level transfer operators, as well as $A_0$.

The number of matrix-free levels $\kappa$ can be adjusted freely, to trade off between the memory cost and the floating point operation count. For our implementation, we set $\kappa=2$, based on our empirical observation that the number of non-zero entries tends to decrease rapidly past this point. Depending on the application, however, higher values of $\kappa$ may be desired, e.g., on GPU architectures, where memory considerations are often a more important consideration for performance.

Note that for simplicity, we omit the details of the matrix-free implementation from our pseudocode in Algorithm \ref{alg:full_mg}.

\section{Numerical examples} \label{sect:examples}
To test the performance of our approach we consider a variety of linear systems $Au=f$ that commonly arise in computational fluid dynamics applications. To assess its performance we apply the multigrid solver both (1) directly and (2) as a left preconditioner for iterative Krylov subspace methods. Either way, we report the iterations required for the relative residual to reach a tolerance of $\frac{\| f-Au \|}{\| f \|} < 10^{-7}$ from an initial guess of zero. Furthermore we record the sizes of the operators at each level of the multigrid hierarchy for a selection of the problems.

The block Jacobi smoother with adaptive damping parameter is used throughout. While we have in our own testing experimented with other smoothers including block Gauss-Seidel, incomplete Cholesky/LU factorisations, in addition to more complicated techniques such as Krylov subspace methods, we use block Jacobi as it is a good choice for many large practical applications as it is extremely simple both to implement and to parallelise.

For all the examples, the number of pre- and post- smoothing steps in the multigrid V-cycle are chosen as $t_\text{pre} = t_\text{post} = 3$. The same number of pre- and post- smoothing steps are used so that the preconditioned Conjugate Gradient (pCG) algorithm can be used for problems involving symmetric positive definite operators. The number of levels is chosen such that the number of aggregates in the coarsest level $L_N$ is exactly equal to $2^d$. As such the number of levels in each multigrid hierarchy for a problem consisting of $|\mathcal{T}_h|$ elements is given by $N=\lceil \, \log_{2^d} |\mathcal{T}_h | \, \rceil$.

\subsection{Elliptic discretisations}
\begin{figure}
    \centering
    \includegraphics[height=50mm]{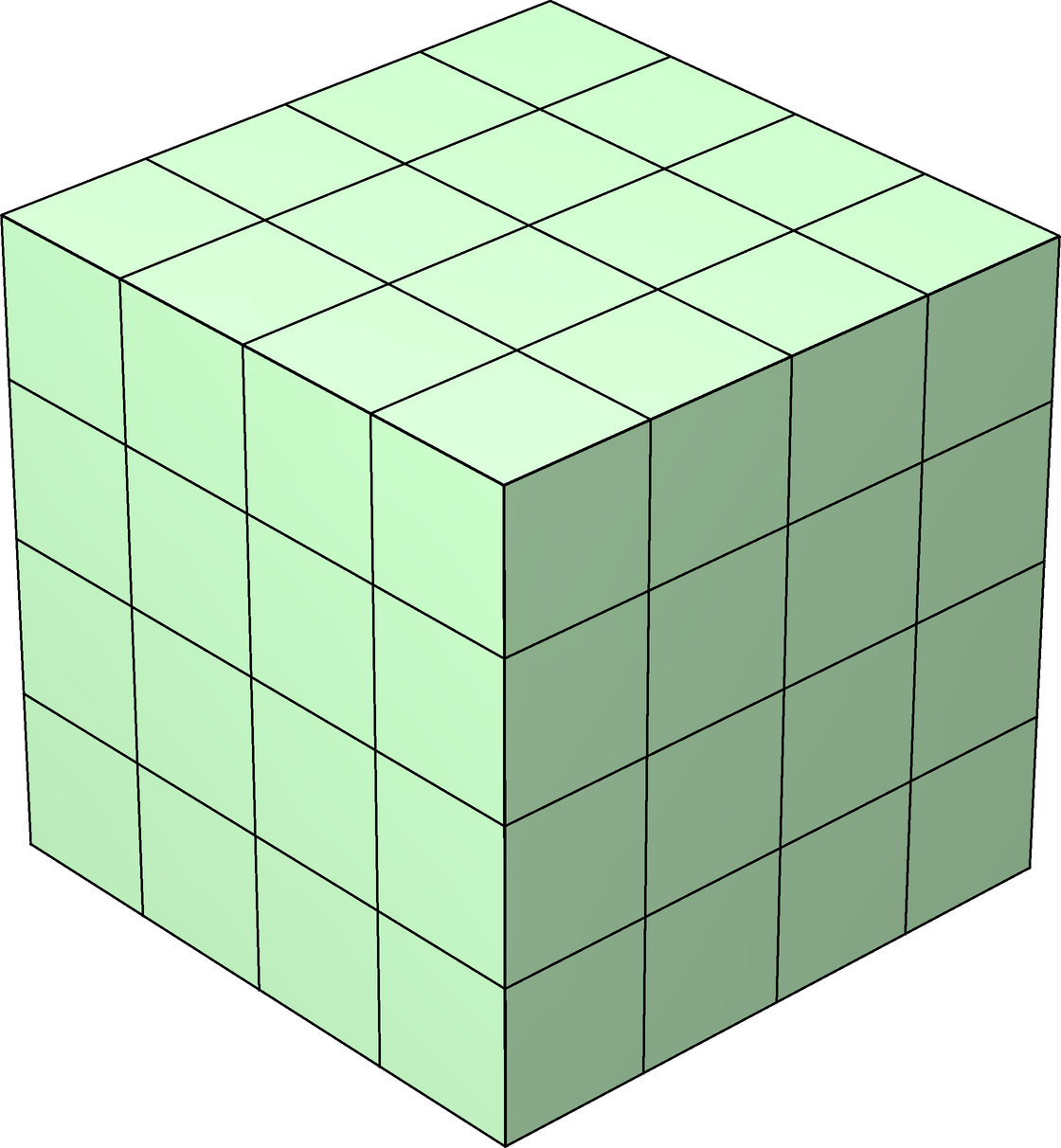}
    \hspace{10mm}
    \includegraphics[height=50mm]{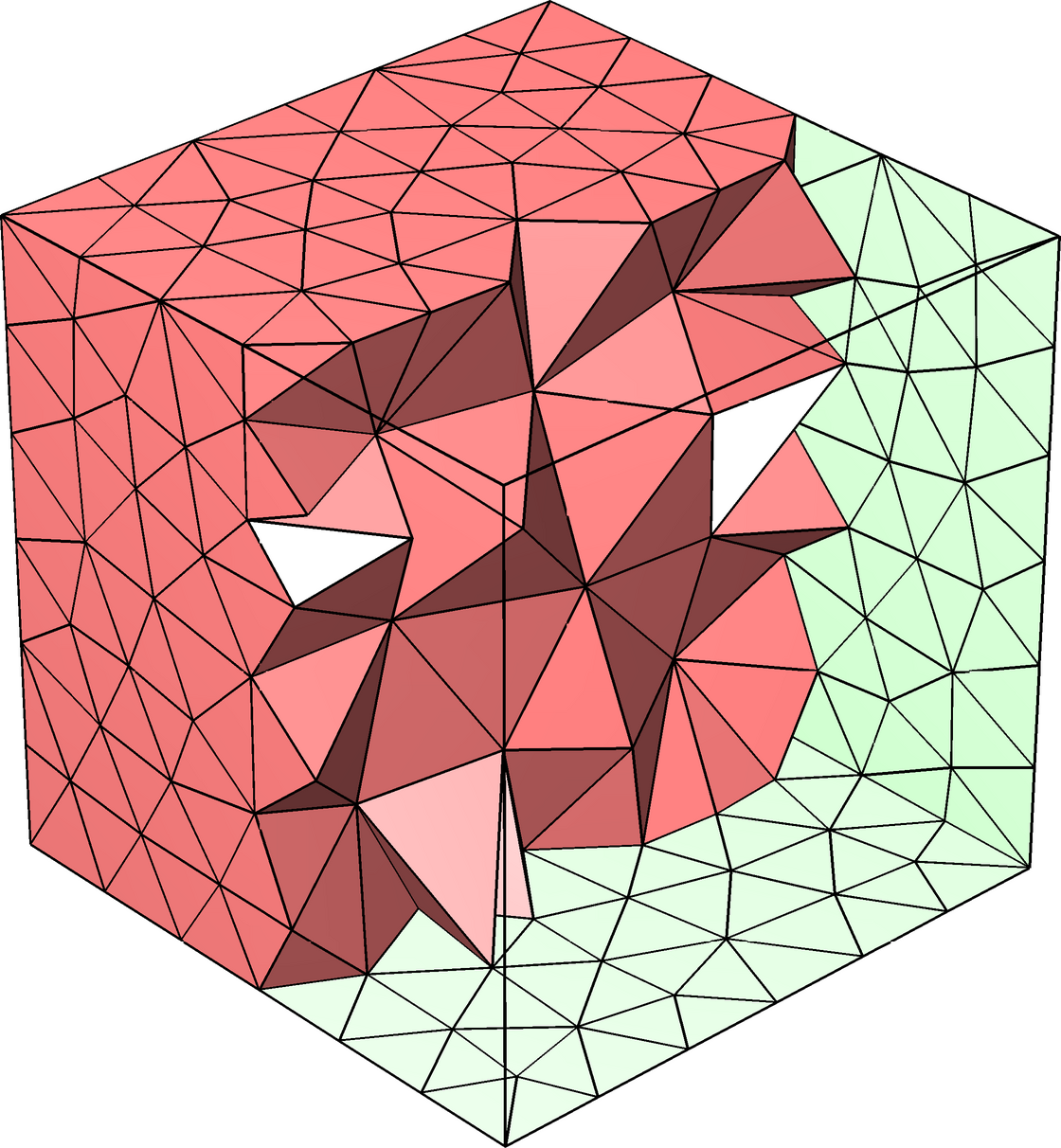}
    \caption{Meshes used for the cube test problems. Left: Structured hexahedral mesh, with refinement level $M=2$. Right: Fully unstructured tetrahedral mesh, at the coarsest level (no refinements). }
    \label{fig:cubemeshes}
\end{figure}

As a first test we first consider Poisson's problem as a model elliptic equation:
\begin{equation}
\label{eq:poisson}
    -\Delta u = f
\end{equation}
on the domain $[-1,1]^3$ in $\mathbb{R}^3$. The domain is discretised using a regular Cartesian mesh consisting of $2^M$ hexahedral elements in each spatial dimension, endowed with outer product polynomials of degree $p=1$. A sample mesh with $M=2$ is shown in on the left of Figure~\ref{fig:cubemeshes}.

In contrast with the scenario for continuous Finite Element methods, there exist multiple popular means of discretising this equation using DG methods. As shown in \cite{arnold2002unified}, all the methods are known to produce consistent and stable Laplacians while differing in their choices of numerical flux. We test our approach using both the Interior Penalty (IP) \cite{arnold1982interior} and minimal dissipation Local Discontinuous Galerkin (LDG) \cite{cockburn1998local} methods. For the interior penalty method the penalty parameter is set to $\mu = \frac{M(p+1)^2}{2}$.

For both discretisations, we choose the right-hand side $f$ such that 
\begin{equation}
    u(x,y,z) = \exp( \sin(\pi x) \sin(\pi y) \sin(\pi z) ) - 1 \label{eqn:poissonexact}
\end{equation}
defines an exact solution for \eqref{eq:poisson} with Dirichlet boundary conditions applied on the left boundary $x=-1$ and Neumann boundary conditions elsewhere. In Figure \ref{fig:poi} below, we report the number of iterations required to converge our approach to the desired relative tolerance of $10^{-7}$, 
applied both directly and as a preconditioner for the Conjugate Gradient algorithm. We also report the sizes of the operators at each level of the multigrid hierarchy for $M=4$.

\begin{figure}[h]
    \centering



    
    \begin{tabular}{ |c||c|c|c|c|  }
     \hline
     $M$ & IP & LDG & IP(pCG) & LDG(pCG) \\
     \hline
     3 & 7 & 7 & 6 & 6 \\
     4 & 7 & 8 & 6 & 7 \\
     5 & 7 & 9 & 6 & 7 \\
     \hline
    \end{tabular}
    \hspace{3mm}
    \begin{tabular}{ |c||r|r|c|  }
     \hline
         $k$ & \multicolumn{1}{|c|}{dof} & \multicolumn{1}{|c|}{nnz} \\
     \hline
     0  & 32,768  & 1,835,008 \\
     1  & 5,206 & 1,706,338 \\
     2  & 882 & 310,538 \\
     3  & 100 & 9,442 \\
     \hline
    \end{tabular}
    \caption{Performance of multigrid solver on structured Poisson example using IP, LDG methods. On the left, we plot the number of iterations for the multigrid method to converge as a solver and as a preconditioner to CG (pCG) for both choices of DG flux, with the corresponding table shown in the centre. On the right, we report the number of unknowns (dof) and non-zero entries (nnz) at each level of the multigrid hierarchy for IP with $M=4$.}
    \label{fig:poi}
\end{figure}

We observe consistent performance of the multigrid construction across both DG fluxes, suggesting that the method is robust to choice of numerical flux. By contrast, it is has been shown in the application of standard geometric multigrid solvers to DG discretisations \cite{fortunato2019efficient, pan2022agglomeration} that differing coarsening strategies for IP and LDG type discretisations must be employed to achieve comparable iteration counts across the board.

Also observed is that the number non-zero entries at the first level of the multigrid hierarchy remains roughly constant despite the rapid decrease in size of the operator. This motivates the use of the matrix-free approach outlined in Section \ref{sect:matrix_free} to circumvent the cost of having to explicitly store the operator in memory.


\subsection{h*- vs h-multigrid}
Next we compare the performance of the h*- and h-multigrid constructions described in Section \ref{sect:hp}, considering the same example problem on $[-1,1]^3 \subset \mathbb{R}^3$ as in the preceding section and using a Cartesian mesh with $2^M$ hexahedral elements in each dimension. For this test we fix $M=4$ and instead vary the polynomial degrees of basis functions used. Furthermore we consider only the minimal dissipation LDG method, for simplicity, as the choice of discretisation.

Figure \ref{fig:poi_hp} shows the number of iterations to reach the desired tolerance using the h*- and h-multigrid solvers both directly and as a preconditioner to the Conjugate Gradient (CG) method. We also report in Figure \ref{fig:poi_hp_nnz} the number of non-zero entries in each level of the multigrid hierarchy for the quadratic $p=2$ case.

\begin{figure}[h]
    \centering



    
    \begin{tabular}{ |c||c|c|c|c|  }
     \hline
     $p$ & h*- & h- & h*-(pCG) & h-(pCG) \\
     \hline
     1 & 8 & 8 & 7 & 6 \\
     2 & 8 & 9 & 6 & 7 \\
     3 & 10 & 9 & 7 & 7 \\
     \hline
    \end{tabular}
    \caption{Performance of h*- and h- multigrid on structured Poisson example using LDG for polynomial degrees $p=1,2,3$. The number of iterations to convergence for each approach as a solver and as a preconditioner to CG are provided in the table.}
    \label{fig:poi_hp}
\end{figure}

\begin{figure}[h]
    \centering
    \begin{minipage}{.4\textwidth}
        \centering
        \caption*{h*-multigrid}
        \begin{tabular}{ |c||r|r|c|  }
         \hline
         $k$ & \multicolumn{1}{|c|}{dof} & \multicolumn{1}{|c|}{nnz} \\
         \hline
         0  & 110,592 & 20,901,888 \\
         1  & 56,029 & 5,411,417 \\
         2  & 9,574 & 5,387,940 \\
         3  & 3,524 & 640,246 \\
         \hline
        \end{tabular}
    \end{minipage}%
    \begin{minipage}{.4\textwidth}
        \centering
        \caption*{h-multigrid}
        \begin{tabular}{ |c||r|r|c|  }
         \hline
         $k$ & \multicolumn{1}{|c|}{dof} & \multicolumn{1}{|c|}{nnz} \\
         \hline
         0  & 110,592 & 20,901,888 \\
         1  & 18,500 & 16,494,040 \\
         2  & 3,104 & 2,342,854 \\
         \hline
        \end{tabular}
    \end{minipage}%
    \caption{Comparison of multigrid hierarchies using h*- and h- multigrid, with quadratic ($p=2$) functions in each element and LDG flux. We record the number of degrees of freedom and non-zero entries at each level $k$ for each case. While h*-multigrid uses one more level, the number of non-zeros past the top level decreases much more rapidly than in the h-multigrid case.}
    \label{fig:poi_hp_nnz}
\end{figure}

The number of iterations for both the h*- and h-multigrid approach appear to remain largely constant as the polynomial degree $p$ increases. We also see that the h*-multigrid approach more rapidly decreases the number of non-zero entries from the top level to the first level, at the cost of introducing one extra level in the hierarchy. In practice, we expect that the best approach may depend on the specific requirements of the user.

\subsection{Varying $n_\text{cut}$} \label{sec:ncut}
Next we examine the effect of varying $n_\text{cut}$ in Equation \eqref{eq:sik}, again by testing our method on Poisson's equation in $\mathbb{R}^3$. The domain $[-1,1]^3$ is discretised using unstructured tetrahedral elements generated using the Gmsh package \cite{geuzaine2009gmsh}. The coarsest mesh (zero refinements) is shown in the right of Figure~\ref{fig:cubemeshes}. We again fix minimal dissipation LDG as our choice of DG flux and the polynomial degree to be $p=1$, while varying the number of refinements nRef of the unstructured mesh. We also again choose the right hand side such that the solution satisfies our model solution (\ref{eqn:poissonexact}), however here we set Dirichlet conditions on the entire domain boundary. As we are keeping the polynomial degree $p=1$ fixed, we are only considering in this example the effect of $n_\text{cut}$ on geometric coarsening. In other words we are only focused on the case where the number of elements in each aggregate on successive levels $k$ and $k+1$ are unequal h-multigrid case in Equation \eqref{eq:sik}.

We report the number of iterations to convergence using the multigrid method only as a preconditioner to CG while varying $n_\text{cut}$. These are shown in Figure \ref{fig:poi_ncut}. We also record the number of degrees of freedom and non-zero entries of each operator at each level of the multigrid hierarchy for the case $\mathrm{nRef}=2$ in Figure \ref{fig:poi_ncut_nnz}.

\begin{figure}[h]
    \centering



    
    \begin{tabular}{ |c||r|r|r|r|r|r|r|r|r|  }
     \hline
     \backslashbox{nRef}{nCut} & $4$ & $5$ & $6$ & $7$ & $8$ \\
     \hline
0 & 7 & 7 &  8 &  8 &  9 \\
1 & 8 & 8 &  9 & 10 & 11 \\
2 & 9 & 9 & 10 & 11 & 11 \\
     \hline
    \end{tabular}
    \caption{Study of multigrid performance varying $n_\text{cut}$ for 3D Poisson's equation using an unstructured tetrahedral mesh. For this example the polynomial order is fixed at $p=1$, and the mesh is refined geometrically nRef times. The number of iterations to convergence for $n_\text{cut}$= 4, 5, ..., 8 at different mesh refinement levels nRef are shown in a table. }
    \label{fig:poi_ncut}
\end{figure}

\begin{figure}[h]
    \centering
    \begin{minipage}{.3\textwidth}
        \centering
        \caption*{$n_\text{cut}=4$}
        \begin{tabular}{ |c||r|r|c|  }
         \hline
         $k$ & \multicolumn{1}{|c|}{dof} & \multicolumn{1}{|c|}{nnz} \\
         \hline
         0 & 214,260 & 11,141,520 \\   
         1 & 53,565 & 20,294,239 \\
         2 & 13,552 & 9,160,578 \\
         3 & 3,403 & 1,437,061 \\
         \hline
        \end{tabular}
    \end{minipage}%
    \begin{minipage}{.3\textwidth}
        \centering
        \caption*{$n_\text{cut}=6$}
        \begin{tabular}{ |c||r|r|c|  }
         \hline
         $k$ & \multicolumn{1}{|c|}{dof} & \multicolumn{1}{|c|}{nnz} \\
         \hline
         0 & 214,260 & 11,141,520 \\   
         1 & 37,099 & 9,829,883 \\
         2 & 6,359 & 2,305,047 \\
         3 & 1,084 & 195,800 \\
         \hline
        \end{tabular}
    \end{minipage}%
    \begin{minipage}{.3\textwidth}
        \centering
        \caption*{$n_\text{cut}=8$}
        \begin{tabular}{ |c||r|r|c|  }
         \hline
         $k$ & \multicolumn{1}{|c|}{dof} & \multicolumn{1}{|c|}{nnz} \\
         \hline
         0 & 214,260 & 11,141,520 \\   
         1 & 28,034 & 5,640,204 \\
         2 & 3,587 & 822,155 \\
         3 & 451 & 42,673 \\
         \hline
        \end{tabular}
    \end{minipage}
    \caption{Study of number of degrees of freedom (dof) and number of non-zero entries (nnz) at each level of multigrid hierarchy for 3D Poisson's equation using an unstructured tetrahedral mesh, while varying $n_\text{cut}$. For these numbers the refinement level is fixed at nRef = 2.}
    \label{fig:poi_ncut_nnz}
\end{figure}

The results show that while the number of iterations needed for convergence is lower for smaller values of nCut, the size of the operators and number of non-zero entries observed in each level of the multigrid hierarchy is higher. Our default choice of $n_\text{cut} = 2^d - d + 1$ aims strikes a balance between iteration count and operator complexity in the multigrid hierarchy. However, depending on user preferences, the value of $n_\text{cut}$ can be varied for desired performance depending on the application.

\subsection{Convection-diffusion comparison}
To benchmark the method against existing multigrid-based alternatives, we consider the stationary convection-diffusion problem
\begin{equation}
    \mu \Delta u + v \cdot \nabla u = 0,
\end{equation}
in which we vary the diffusion coefficient $\mu$ and the velocity $v = (1,2,3)/\sqrt{14}$ is a specified velocity. In addition to our method, we consider the use of block Jacobi, classical AMG, and smoothed aggregation AMG preconditioners for the GMRES iterative method to solve this problem. For the two AMG implementations we use the AlgebraicMultigrid package \cite{amg_julia} available in the Julia programming language.

For the DG discretisation, we use the minimal dissipation LDG method for the diffusion term and the Godunov flux for the convection term. The problem is discretised on the domain $[-1,1]^3$ again using unstructured tetrahedral meshes from Gmsh. For this test, we vary only the number of mesh refinements and use linear ($p=1$) elements. With this setup we report the number of preconditioned GMRES iterations using each of the aforementioned preconditioners for three different P\'eclet numbers $\mathrm{Pe}=\|v\|L/\mu = 2/\mu$ shown in Figure \ref{fig:conv_diff_0}. The GMRES restart parameter is always set to 1000.

\begin{figure}[h]
    \centering



    
    \begin{tabular}{ |c|r||r|r|r|r|r|r|r|  }
     \hline
     nRef & dof & Our method & block Jacobi & AMG(classical) & AMG(SA) \\
     \hline
     \multicolumn{6}{|l|}{\hspace{4.2cm}Pe$~= 0$ (pure viscosity)} \\
     \hline
     0 & 4,188   & 9 &  184 &  38 &  36 \\
     1 & 29,180  & 10 &  430 &  64 &  58 \\
     2 & 214,260 & 11 & 1462 & 121 &  91 \\
     3 & 1,660,556 & 13 & 3838 & 309 & 185 \\
     \hline
     \multicolumn{6}{|l|}{\hspace{4.2cm}Pe$~= 100$} \\
     \hline
     0 & 4,188 &  9 &  76  & 61 &  66 \\
     1 & 29,180 &  10 & 181  & 57  & 74 \\
     2 & 214,260 & 12 & 451 &  75 & 100 \\
     3 & 1,660,556 & 16 & 892 & 173 & 172 \\
     \hline
     \multicolumn{6}{|l|}{\hspace{4.2cm}Pe$~= 1,000$} \\
     \hline
     0 & 4,188 & 10 &  46 &  -  &  - \\
     1 & 29,180 & 12 &  83 &  - & 2473 \\
     2 & 214,260 & 12 & 279 & 719 &  507 \\
     3 & 1,660,556 &  17 & 589 & 367 &  325 \\
     \hline
     \multicolumn{6}{|l|}{\hspace{4.2cm}Pe$~ = \infty$ (pure convection)} \\
     \hline
     0 & 4,188 & 9 & 37 & - & - \\
     1 & 29,180 & 13 & 63 & - & - \\
     2 & 214,260 & 17 & 118 & - & - \\
     3 & 1,660,556 &  25 & 416 & - & - \\
     \hline
    \end{tabular}
    \caption{Comparison of number of GMRES iterations under h-refinement using different preconditioners for convection-diffusion problem with varying diffusion coefficients. For this problem the domain is discretised with an unstructured tetrahedral mesh. A dash symbol (-) is used for the cases where GMRES did not converge within 5000 iterations.}
    \label{fig:conv_diff_0}
\end{figure}

The results suggest for our solver to be applicable across different convection-diffusion regimes. In particular we note that these results were obtained using a block Jacobi smoother across the board and that in practice more advanced alternative smoothers such as incomplete LU factorisations as applied in \cite{persson2008newton} for convection dominated cases can be used to obtain further improvements in solver performance.

In contrast, for the classical and smoothed aggregation AMG solvers considered here, the solvers appear to perform well only in regimes where the problem becomes diffusion-dominated (i.e., small P\'eclet number). Moreover, for these alternative approaches, the number of iterations to convergence appears to grow much more severely as the mesh is refined spatially, even in the diffusion-dominated regime, compared to our method. 

\section{Conclusion} \label{sect:conclusion}
In this work we have introduced a geometric multigrid solver for Discontinuous \\ Galerkin discretisations using ideas from adaptive smoothed aggregation AMG. The solver is simple to construct and can be employed in a matrix-free fashion enabling it to be deployed on distributed computing architectures. Numerical tests validate that the method is applicable to both first- and second-order partial differential equations and can be used both as an iterative solver and as a preconditioner. Remarkably, no individualized tuning is required for different choices of numerical flux, in contrast to classical geometric multigrid techniques.

We will provide open-source code for our methodology at 
\[
\texttt{github.com/yllpan/adaptiveMG}
\]
We plan to maintain the code hosted at the link, and any subsequent work extending the solver will be reflected in the repository. In future work we intend to benchmark our method against existing DG solver packages for a wide range of applications including, but not limited to, computational fluid dynamics. Moreover, we plan to investigate a framework for the analysis of the method with a view toward justifying more carefully and improving upon some of the parameter choices adopted here.

\section*{Acknowledgments}
This material is based on work supported 
by the Applied Mathematics
Program of the US Department of Energy (DOE) Office of Advanced Scientific
Computing Research under contract number DE-AC02-05CH11231 
and by the U.S. Department of Energy, Office of Science, Accelerated Research in Quantum Computing Centers, Quantum Utility through Advanced Computational Quantum Algorithms, grant no. DE-SC0025572. 
M.L. was also partially supported by
a Sloan Research Fellowship, and Y.P. and P.P. were  also partially supported by the National Science Foundation under Grant DMS-2309596.

\newpage
\bibliographystyle{plain}
\bibliography{references}

\end{document}